\documentclass[12pt]{article}

\usepackage{amsmath}
\usepackage{amssymb}
\usepackage{amsthm}
\usepackage[all]{xy}
\usepackage{longtable}

\usepackage[margin=1in]{geometry}

\theoremstyle{plain}
\newtheorem{thm}{Theorem}[section]
\newtheorem{cor}[thm]{Corollary}
\newtheorem{lem}[thm]{Lemma}
\newtheorem{prop}[thm]{Proposition}
\newtheorem*{rmk}{Remarks}

\newcommand{\Q}{\mathcal{Q}}
\newcommand{\q}{\mathbb{Q}}
\newcommand{\C}{\mathfrak{C}}
\newcommand{\mcA}{\mathcal{A}}
\newcommand{\mcL}{\mathcal{L}}
\newcommand{\M}{\mathfrak{M}}
\newcommand{\Cl}{\mathcal{C}l}
\newcommand{\m}{\mathfrak{m}}
\newcommand{\n}{\mathfrak{n}}
\newcommand{\f}{\mathfrak{f}}
\newcommand{\mcF}{\mathcal{F}}
\newcommand{\mcB}{\mathcal{B}}

\newcommand{\bbz}{\mathbb{Z}}
\newcommand{\bbn}{\mathbb{N}}
\DeclareMathOperator{\Ass}{Ass} \DeclareMathOperator{\Supp}{Supp}
\DeclareMathOperator{\I}{I}
\DeclareMathOperator{\cok}{coker}\DeclareMathOperator{\image}{Im}
\DeclareMathOperator{\rank}{rank} \DeclareMathOperator{\glb}{glb}
\DeclareMathOperator{\Gal}{Gal}
\DeclareMathOperator{\depth}{depth}
\DeclareMathOperator{\Spec}{Spec} \DeclareMathOperator{\End}{End}
\DeclareMathOperator{\syz}{syz}
\newcommand{\ses}[3]{\ensuremath{0 \rightarrow #1 \rightarrow #2 \rightarrow #3 \rightarrow 0}}
\newcommand{\D}{\ensuremath{x^2y-y^{n-1}}}
\newcommand{\A}{\ensuremath{x^2-y^{n+1}}}
\newcommand{\Esix}{\ensuremath{x^3-y^4}}
\newcommand{\Eseven}{\ensuremath{x^3-xy^3}}
\newcommand{\Eeight}{\ensuremath{x^3-y^5}}

\newcommand{\R}[1]{\ensuremath{k[[x,y]]/(#1)}}

\everymath={\displaystyle}

\begin{document}

\title{A Krull-Schmidt Theorem for One-dimensional Rings of Finite Cohen-Macaulay Type}
\author{Nicholas R. Baeth\footnote{This work consists of research done as part of the author's Ph.D. thesis, conducted at the University of Nebraska--Lincoln under the direction of Roger Wiegand.} \ \\ Department
of Mathematics \ \\ University of Nebraska--Lincoln \ \\ Lincoln,
NE 68588-0130\ \\ nbaeth@math.unl.edu}
\date{October 28, 2004}

\maketitle

\section{Introduction}

Let $(R,\m,k)$ be a one-dimensional Cohen-Macaulay local ring
(where $\m$ is the maximal ideal and $k=R/\m$). We assume
throughout that $R$ is equicharacteristic, equivalently, $R$
contains a field. Recall that $R$ has finite Cohen-Macaulay type
provided there are, up to isomorphism, only finitely many
indecomposable maximal Cohen-Macaulay modules (finitely generated
torsion-free modules in this setting). If $R$ has finite
Cohen-Macaulay type, the monoid $\C(R)$ of isomorphism classes of
maximal Cohen-Macaulay (MCM) modules (together with $[0]$) is
isomorphic to a submonoid of some $\bbn^t$ defined by a finite
family of homogeneous linear equations with integer coefficients.
That is, there is an $s\times t$ matrix $\mcA$ such that
$\C(R)\cong \bbn^t \cap \ker(\mcA)$. (Here $\bbn$ denotes the
monoid of non-negative integers.)

For each ring $(R,\m,k)$ of finite Cohen-Macaulay type, with $k$
perfect and of characteristic different from $2,3$ and $5$, we
determine exactly the defining equations for the monoid $\C(R)$.
From these defining equations we are able to determine exactly
when $\C(R)$ is free, that is, direct-sum decompositions of MCM
$R$-modules have the Krull-Schmidt uniqueness property. Further,
we determine which rings have the weaker property that any two
representations of a MCM module as a direct sum of indecomposables
have the same number of indecomposable summands.

If $R$ is complete, then $\C(R)$ is the free monoid on the set of
isomorphism classes of indecomposable MCM modules. In order to
describe the monoid $\C(R)$ in the non-complete case we need a
detailed description of the indecomposable MCM $\hat R$-modules,
together with information on their ranks at the various minimal
prime ideals of $\hat R$. When $k$ is algebraically closed, we can
glean this information from the Auslander-Reiten quivers for $\hat
R$, which are worked out in detail in Yoshino's book \cite{Yos90}.
To complete our study of $\C(R)$ in the incomplete case and the
case where $k$ is perfect but not algebraically closed, we analyze
the maps $\C(R) \to \C(\hat R)$ and $\C(R) \to \C(S)$, where
$(S,{\n}, \ell)$ is a flat local extension of $R$ and $\ell/k$ is
an algebraic extension.

\section{The Hierarchy of Complete Rings of FCMT}

In this section we describe all complete one-dimensional
equicharacteristic local Cohen-Macaulay rings of finite
Cohen-Macaulay type (FCMT). We recall the classification, given in
\cite{RW91}. We have changed the names of the rings given in
\cite{RW91} in order to match the more common labels given in the
literature (e.g. \cite{GK85}, \cite{Yos90}).

\begin{thm}\label{classification} Let $(R,\m,k)$ be a complete one-dimensional
equicharacteristic Cohen-Macaulay local ring. Further assume that
$k=R/\m$ is perfect of characteristic not $2$, $3$, or $5$. Then
$R$ has finite Cohen-Macaulay type if and only if $R$ is
isomorphic to:

\begin{enumerate}
\item One of the hypersurfaces $k[[x,y]]/(f)$ listed in Table
\ref{hypersurfaces} or

\item \label{non-hypersurfaces} $\End_S(\n)$ where $(S, \n)$ is
one of the rings listed in Table \ref{hypersurfaces} not of type
$(A_1)$.
\end{enumerate}

\begin{center}
\begin{longtable}{|llrc|}\caption{One-dimensional Hypersurfaces of FCMT}\label{hypersurfaces}\\
\hline \endhead \hline \endfoot
type & $R$ & & $[K:k]$ \\
\hline
$(A_n)$ & \R \A & ($n \geq 0$) & $1$\\
$(D_n)$ & \R \D & ($n \geq 4$) & $1$\\
$(E_6)$ & \R \Esix & & $1$\\
$(E_7)$ & \R \Eseven & & $1$\\
$(E_8)$ & \R \Eeight & & $1$\\
\hline
$(A2_n)$ & $k[[T,\xi T^{n+1}]]$ & ($n \geq 1$)  & $2$\\
$(D2_n)$ & $k[[(T,U),(\xi T^n,U), (0,U^2)]]$ & ($n \geq 1$)  & $2$\\
$(D3)$ & $k[[T, \xi T]]$ & &  $3$\\
\hline
\end{longtable}
\end{center}

\end{thm}

The notation in the table deserves some explanation. We denote the
integral closure of $R$ in its total quotient ring by $\bar R$ and
let $K$ be a residue field of $\bar R$ with maximal degree $[K:k]$
over $k$. In the classification of the hypersurfaces (we will see
shortly that the ring of type ($D2_n$) is indeed a hypersurface)
of FCMT there is a natural dichotomy. When $K=k$, $R$ is a ring of
type ($A_n$), ($D_n$), ($E_6$), ($E_7$), or ($E_8$). The remaining
cases occur when $[K:k]>1$. When $[K:k]=2$ we have the cases
($A2_n$) and ($D2_n$) with $\xi \in K-k$. The isomorphism class of
$R$ is independent of the choice of $\xi\in K$ but varies with
choices of $K$. When $[K:k]=3$ we have the case ($D3$). Again, the
choice of $\xi\in K-k$ does not affect the isomorphism class of
$R$.

We now give defining equations for these three families of rings
as we will need them in the following section. First let $R$ be a
ring of type ($A2_n$). Letting $x=\xi T^{n+1}$ and $y=T$, we have
$\xi=\frac{x}{y^{n+1}}$. Since the isomorphism class of $R$ is not
dependent on which $\xi$ we choose we can take $\xi \in K-k$ such
that $\xi^2 \in k$. Now the minimal polynomial for $\xi$ over $k$
has the form $X^2-\xi^2$. Thus
$\left(\frac{x}{y^{n+1}}\right)^2-\xi^2=0$ and hence
$x^2-\xi^2y^{2n+2}=0$. Then $$R \cong
k[[x,y]]/(x^2-\xi^2y^{2n+2}).$$

Now suppose $R$ is a ring of type ($D2_n$). Again we choose
$\xi\in K$ such that $\xi^2 \in k$. Letting $x=(\xi T^n, U)$,
$y=(T,U)$ and $z=(0,U^2)$, we can write $R$ as
$$R \cong \frac{k[[x,y,z]]}{(\xi^2y^{2n}-x^2-\xi^2z^n+z, z^n-y^{2n-2}z,yz-xz)}.$$

Note that since $z^n=y^{2n-2}z$ and $\xi^2y^{2n}-x^2=\xi^2z^n-z$
we have that $$z=\frac{\xi^2y^{2n}-x^2}{\xi^2y^{2n-2}-1}.$$
Therefore
$$y\left(\frac{\xi^2y^{2n}-x^2}{\xi^2y^{2n-2}-1}\right)-x\left(\frac{\xi^2y^{2n}-x^2}{\xi^2y^{2n-2}-1}\right)=0$$
and so $$R \cong \frac{k[[x,y]]}{(x-y)(x^2-\xi^2y^{2n})}.$$

Finally, suppose that $R$ is a ring of type ($D3$). Let
$X^3+aX^2+bX+c$ be the minimal polynomial for $\xi$ over $k$.
Letting $x=T$ and $y=\xi T$, we have $\xi=\frac{y}{x}$ and hence
$\left(\frac{y}{x}\right)^3+a\left(\frac{y}{x}\right)^2+b\left(\frac{y}{x}\right)+c=0$.
Therefore $y^3+ay^2x+byx^2+cx^3=0$, and
$$R \cong k[[x,y]]/(y^3+ay^2x+byx^2+cx^3).$$

We adopt the following notation for the non-hypersurface rings of
part \ref{non-hypersurfaces} in Theorem \ref{classification}: we
say that $(R',\m',k):=\End_R(\m)$ is a ring of type ($D'_n$)
(resp. ($E'_6$), ($E'_7$), ($E'_8$), ($D2'_n$), ($D3'$)) if
$(R,\m,k)$ is a ring of type ($D_n$) (resp. ($E_6$), ($E_7$),
($E_8$), ($D2_n$), ($D3$)) listed in Table \ref{hypersurfaces}. We
note that if $R$ is a ring of type ($A_n$) with $n \geq 2$ then
$R'=\End_R(\m)$ is a ring of type ($A_{n-2}$) and that if $R$ is a
ring of type ($A2_n$) with $n \geq 1$ then $R'=\End_R(\m)$ is a
ring of type ($A2_{n-1}$).

When $k$ is an algebraically closed field of characteristic zero,
Theorem \ref{classification} gives the classification in
\cite{GK85}. We now express the rings from \cite{GK85} as flat
local extensions of rings in Table \ref{hypersurfaces}. Let $R$ be
a ring of type ($A2_n$) and let $K$ be as in the paragraph after
Theorem \ref{classification}. Then
\begin{equation}\label{A2}S=R \otimes_kK \cong K[[x,y]]/((x-\xi y^{n+1})(x+\xi
y^{n+1})),\end{equation} which is a ring of type ($A_{2n+1}$).
Similarly, if $R$ is a ring of type ($D2_n$), then
\begin{equation}\label{D2}S=R\otimes_kK \cong K[[x,y]]/((x-y)(x-\xi y^n)(x+\xi
y^n)),\end{equation} which is a ring of type ($D_{2n+2}$). Now if
$R$ is a ring of type ($D3$) and $L$ is the Galois closure of
$K/k$ we see that
\begin{equation}\label{D3}S=R \otimes_kL \cong L[[x,y]]/(y^3+ay^2x+byx^2+cx^3),\end{equation} which is a
ring of type ($D_4$) since $y^3+ay^2x+byx^2+cx^3$ splits into
linear terms over $L$.

Now if $R \cong k[[x,y]]/(f)$ is a ring of type ($A_n$), ($D_n$),
($E_6$), ($E_7$), or ($E_8$) with $k$ perfect and not
algebraically closed we note that
\begin{equation}\label{algclose}k[[x,y]]/(f) \longrightarrow \bar
k[[x,y]]/(f),\end{equation} where $\bar k$ is the algebraic
closure of $k$, is a flat local extension of rings. The flatness
of this map follows from the fact that $\bar k[x,y]$ is faithfully
flat over $k[x,y]$ and from \cite[Thm. 22.4]{Mat00}.

\section{Indecomposable MCM Modules}

In this section we classify all of the indecomposable maximal
Cohen-Macaulay modules over the rings listed in Theorem
\ref{classification}. When $k$ is algebraically closed one can
find the classification in \cite{GK85} and \cite[Ch. 9]{Yos90}.
When $k$ is not algebraically closed we consider faithfully flat
extensions $R \rightarrow (S,n,\ell)$ where $S$ is as in
(\ref{A2}), (\ref{D2}), (\ref{D3}), or (\ref{algclose}).

Let $\C(T)$ denote the set of isomorphism classes of MCM modules
over a ring $T$. We consider the map on modules $M \mapsto M
\otimes_R S$. Since $S$ is faithfully flat as an $R$-module, this
map is one-to-one up to isomorphism (c.f \cite[2.5.8]{EGA4}).
Moreover, taking $N=S$ in \cite[Thm. 2.17]{BH98}, we see that this
map takes MCM $R$-modules to MCM $S$-modules. Therefore, we have
an injection $\C(R) \hookrightarrow \C(S)$. Since we know the MCM
$S$-modules it is enough to determine which MCM $S$-modules are
extended ($M \cong N \otimes_R S$ for some $R$-module $N$) in
order to classify all MCM $R$-modules. Moreover, the following
result allows us to conclude that extended $S$-modules with no
extended proper direct summands are extended from indecomposable
$R$-modules. We recall that the Krull-Schmidt property holds for
complete local rings, \cite{A50}, and direct sum cancellation
holds over all local rings, \cite{E73}. We write $M \mid N$ to
indicate that $M$ is isomorphic to a direct summand of $N$.

\begin{lem}\label{divhom}
Let $R \rightarrow S$ be a faithfully flat extension of Noetherian
rings. Suppose that the Krull-Schmidt theorem holds for finitely
generated modules over $R$ and direct sum cancellation of finitely
generated modules holds over $S$. Then, given finitely generated
$R$-modules $M$ and $N$,
$$M\mid N \text{ as $R$-modules if and only if }(S
\otimes_RM)\ \big |\ (S\otimes_RN) \text{ as $S$-modules}.$$
\end{lem}

\begin{proof}
It is clear that $M\mid N$ implies $(S \otimes_RM)\ \big |\
(S\otimes_RN)$. Now suppose that $(S\otimes_RM)\ \big |\
(S\otimes_RN)$. From \cite[Lem. 1.2]{RW98} we know that $M\mid
N^r$ for some $r>1$. We induct on the length of the direct-sum
decomposition of $M$. If $M$ is indecomposable, $M\mid N$ by the
Krull-Schmidt theorem. Otherwise, write $M =M_1 \oplus M_2$ where
$M_1$ is indecomposable. Then $M_1\mid N$, say $N\cong M_1 \oplus
N_1$. Now $$[(S \otimes_R M_1) \oplus (S \otimes_R M_2)] \ \big |\
[(S \otimes_R M_1) \oplus (S \otimes_R N_2)],$$ and by
cancellation over $S$, $(S \otimes_R M_2) \ \big |\  (S \otimes_R
N_2)$. By the induction hypothesis $M_2 \mid N_2$, and hence $M
\mid N$.
\end{proof}

We first consider $R$ to be a ring of type ($A_n$), ($D_n$),
($E_6$), ($E_7$), or ($E_8$). If, in addition, $R$ has
algebraically closed residue field $k$, we turn to \cite[Ch.
9]{Yos90} for a complete description (by way of matrix
factorizations) of each of the indecomposable MCM $R$-modules.

If $R$ is a ring of type ($A_n$), ($D_n$), ($E_6$), ($E_7$), or
($E_8$) and $k$ is not algebraically closed we consider the flat,
local extension (\ref{algclose}) above and the following lemma.

\begin{lem}\label{close} Let $R=k[[x,y]]/(f)$ and $S=\bar k[[x,y]]/(f)$ be
rings of (the same) type ($A_n$), ($D_n$), ($E_6$), ($E_7$), or
($E_8$). Then there is a one-to-one correspondence between the
indecomposable MCM $R$-modules and the indecomposable MCM
$S$-modules given by $_RM \mapsto M \otimes_R S$.
\end{lem}

\begin{proof} The indecomposable MCM $S$-modules are explicitly
calculated in \cite[Ch. 9]{Yos90} as the cokernels of certain
matrices. The entries of these matrices are all monomials in $x$
and $y$ (where $\m=(x,y)$ is the maximal ideal of $R$) with
coefficients in $k$. (We note that although the matrix
factorizations of Yoshino involve $i=\sqrt{-1}$, we can avoid this
by using the equations in Table \ref{hypersurfaces} rather than
those in \cite{Yos90}. The alternate matrix factorizations are
worked out in \cite{NB05}.) Thus every indecomposable MCM
$S$-module $N$ is extended from a finitely generated $R$-module
$M$. Obviously $M$ is indecomposable, and $M$ is a MCM module by
\cite[Prop. 1.2.16, Thm. A.11]{BH98}.\end{proof}

For the remainder of this section $R$ is a ring of type ($A2_n$),
($D2_n$), or ($D3$). Furthermore, if $R$ is as in (\ref{A2}) or
(\ref{D2}) we put $S=R \otimes_kK$ and if $R$ is as in (\ref{D3})
we put $S=R\otimes_kL$.

We note that $R \rightarrow S$ is a flat local extension of rings
and by Lemma \ref{close} we know exactly the MCM $S$-modules. We
will now determine which of the MCM $S$-modules are extended from
MCM $R$-modules and thus determine the MCM $R$-modules.

We begin by stating the following proposition, whose proof
consists of the essential details in the proof of the
Krull-Schmidt theorem (c.f. \cite[Sec. 5.4]{P82}). Recall that a
ring $E$ is ``local'' in the non-commutative sense provided that
$E/J(E)$ (where $J(E)$ is the Jacobson radial of $E$) is a
division ring, equivalently, has a unique maximal left ideal.

\begin{prop} Let $(R,\m)$ be a local ring and suppose $X_1 \oplus
\cdots X_s \cong Y_1 \oplus \cdots \oplus Y_t$, where $X_i$ and
$Y_j$ are indecomposable finitely generated $R$-modules. Suppose
$\End_R(X_1)$ has a unique maximal left ideal. Then $X_1 \cong
Y_j$ for some $j$. \end{prop}

Since direct sum cancellation holds over local rings \cite{E73} we
have:

\begin{cor}\label{partKS} Let $(R,\m)$ be a local ring and $M$ a finitely generated
$R$-module. If $\End_R(X)$ has a unique maximal left ideal and $X$
occurs with multiplicity $\mu$ in some direct sum decomposition of
$M$ then $X$ occurs with multiplicity $\mu$ in every direct sum
decomposition of $M$. \end{cor}

We note in particular that if $C$ is a cyclic module over a local
ring $R$ then $\End_R(C)$ is local. We define the
\emph{multiplicity of $C$ in $M$}, $\mu_M(C)$, to be the number of
copies of $C$ that occur in a direct sum decomposition of $M$. By
Corollary \ref{partKS}, $\mu_M(C)$ is well-defined.

To simplify the notation, we write $L$ in place of $K$ when $R$ is
a ring of type ($A2_n$) or ($D2_n$). Then, when $R$ is a ring of
type ($A2_n$), ($D2_n$) or ($D3$), we have that $R \rightarrow
S:=R\otimes_kL$ is a Galois extension of rings with Galois group
$G:=\Gal_R(S)=\Gal_k(L)$ (c.f. \cite[Prop. 2.4]{RS82}). That is,
(1) $R$ is the $G$-invariant subring of $S$ and (2) for all
subgroups $H$ of $G$ and all $H$-stable ideals $I$ of $S$ with $I
\not= S$, $H$ acts faithfully on $S/I$. Recall the following fact
about Galois extensions. A proof may be found in \cite[Sec.
2.2]{B83}.

\begin{prop}\label{transitive}
Let $A\subseteq B$ be a finite Galois extension of rings with
Galois group $G$. Let $Q$ be a prime ideal of $A$ and let
$\mathcal{P}=\{P_1, \ldots , P_s\}$ be the set of primes of $B$
lying over $Q$. Then $G$ acts transitively on $\mathcal{P}$.
\end{prop}

We say that an action of the Galois group $G$ on an $S$-module $M$
is \emph{semi-linear} if for each $g \in G$,
$$g(sm)=g(s)g(m)$$ for all $s \in S$ and $m \in M$. We know that
the $S$-module $M$ is extended if and only if G acts semi-linearly
on $M$. Indeed, if G acts semi-linearly on $M$, then $M \cong S
\otimes_R M^G$ where $M^G$ is the fixed module (cf. \cite[Prop.
2.5]{RS82}). Conversely, if $M=S \otimes_RN$ for some $R$-module
$N$, then $G$ acts semi-linearly on $M$ via $g(s \otimes
n)=g(s)\otimes n$. The following proposition allows us to
determine the cyclic MCM $R$-modules.

\begin{prop}\label{galois}
Let $A \subset B$ be a Galois extension of rings with Galois group
$G$ ($A=B^G$). Let $Q$ be a prime of $A$ and let $\mathcal
P=\{P_1, \ldots , P_s\}$ be the set of primes of $B$ lying over
$Q$.

\begin{enumerate}
\item Suppose $s=2$ and $G=\langle g \rangle$ where $g$ has order
two. Let $P$ be a prime of $A$ such that $g(P)=P$. Then the
$B$-modules $B/P_1 \oplus B/P_2$ and $B/(P_1 \cap P) \oplus B/(P_2
\cap P)$ are extended from $A$-modules. Moreover, if $M$ is an
extended $B$-module, $\mu_M(B/P_1)=\mu_M(B/P_2)$ and
$\mu_M(B/(P_1\cap P))=\mu_M(B/(P_2\cap P))$.

\item Let $s=3$. Then the $B$-modules $B/P_1 \oplus B/P_2 \oplus
B/P_3$ and $B/(P_1 \cap P_2) \oplus B/(P_1 \cap P_3) \oplus B/(P_2
\cap P_3)$ are extended from $A$-modules. Moreover, if $M$ is an
extended $B$-module, then $\mu_M(B/P_1)=\mu_M(B/P_2)=\mu_M(B/P_3)$
and $\mu_M(B/(P_1 \cap P_2))=\mu_M(B/(P_1 \cap P_3))=\mu_M(B/(P_2
\cap P_3))$.
\end{enumerate}
\end{prop}

\begin{proof}\begin{enumerate} \item We have $g(P_1)=P_2$ and $g(P_2)=P_1$ by
Proposition \ref{transitive}. Therefore $g$ has a natural action
on the $B$-module $B/P_1 \oplus B/P_2$ sending $(b_1 + P_1, b_2 +
P_2)$ to $(g(b_2)+P_1, g(b_1)+P_2)$. One checks easily that this
action is semi-linear. Thus $B/P_1 \oplus B/P_2$ is an extended
module.

If $M$ is any extended $B$-module then $G=\langle g \rangle$ acts
semi-linearly on $M$. Suppose that $_BM$ decomposes into a direct
sum of $B$-modules $N \oplus N'$. It is easily checked that $g(N)$
and $g(N')$ are $B$-submodules of $M$ and that $g(N)\oplus
g(N')=M$. If $N \cong B/P_1$ the composition of the three maps in
the following diagram gives a $B$-module isomorphism $g(N) \cong
B/P_2$. $$g(N) \xrightarrow{g^{-1}} N \xrightarrow{\cong} B/P_1
\xrightarrow{g} B/P_2$$ Therefore, if $M$ has a direct summand
isomorphic $B/P_1$, $M$ also has a direct summand isomorphic to
$B/P_2$. A symmetric argument shows that if $M$ has a direct
summand isomorphic to $B/P_2$, $M$ also has a direct summand
isomorphic to $B/P_1$. Thus $\mu_M(B/P_1) = \mu_M(B/P_2)$. The
other assertions are proved similarly.

\item Now suppose that $s=3$. For each $i\in \{1,2,3\}$ and each
$g \in G$ there is a unique index $g(i) \in \{1,2,3\}$ such that
$g(P_i)=P_{g(i)}$. Then $G$ has a natural action on the $B$-module
$B/P_1 \oplus B/P_2 \oplus B/P_3$ defined by $$g(b_1+P_1, b_2+P_2,
b_3+P_3) \mapsto
(g(b_{g(1)})+P_1,g(b_{g(2)})+P_2,g(b_{g(3)})+P_3).$$ This action
is semi-linear, and hence  the $B$-module $B/P_1 \oplus B/P_2
\oplus B/P_3$ is extended.

Let $M$ be any extended $B$-module and let $g \in G$. The argument
in part one shows that if $N$ is a $B$-direct summand of $M$, so
is $g(N)$. Using Proposition \ref{transitive}, choose $g,h \in G$
be such that $g(P_1)=P_2$ and $h(P_2)=P_3$. If $N \cong B/P_1$,
the diagrams
$$g(N) \xrightarrow{g^{-1}} N \xrightarrow{\cong} B/P_1
\xrightarrow{g} B/P_2$$ and $$hg(N) \xrightarrow{(hg)^{-1}} N
\xrightarrow{\cong} B/P_1 \xrightarrow{hg} B/P_3$$ give
isomorphisms $g(N) \cong B/P_2$ and $hg(N) \cong B/P_3$. Thus if
$M$ has a direct summand isomorphic to $B/P_1$, $M$ also has
direct summands isomorphic to $B/P_2$ and $B/P_3$. Using symmetry
we conclude that $\mu_M(B/P_1)=\mu_M(B/P_2)=\mu_M(B/P_3)$. The
other assertions are proved similarly.

\end{enumerate}\end{proof}

In what follows we exhibit matrix factorizations over $S$ which
are also matrix factorizations over $R$. Thus we show that the MCM
$S$-modules corresponding to the cokernels of these matrices are
extended from MCM $R$-modules. We first make the following
remarks, most of which may be found in Chapter 7 of \cite{Yos90}.

\begin{rmk} Let $B$ be a regular local ring and let $f$ be a
non-zero non-unit of $B$. Then let $A=B/(f)$, a hypersurface.
Recall that a matrix factorization of $f$ is a pair of square ($n
\times n$) matrices $(\phi, \psi)$ with entries in $B$ such that
$\phi\cdot\psi=f\cdot 1_{B^n}$ and $\psi\cdot\phi=f\cdot 1_{B^n}$.
We now recall some basic facts about matrix factorizations.

\begin{enumerate}

\item We will only consider reduced matrix factorizations, that
is, matrix factorizations given by matrices whose entries are not
units. The cokernel of such a matrix has no non-zero free direct
summand.

\item If $(\phi, \psi)$ is a reduced matrix factorization over a
ring $R$, $\cok(\psi) \cong \syz^1_R(\cok(\phi))$.

\item If $(\phi, \psi)$ is a reduced matrix factorization of $f$,
then $B^n \xrightarrow{\phi} B^n \rightarrow \cok(\phi)$ is a
minimal free presentation for $\cok(\phi)$. Thus, if $(\phi',
\psi')$ is another reduced matrix factorization of $f$, then
$\I_1(\cok(\phi))=\I_1(\cok(\phi')$ and
$\I_2(\cok(\phi))=\I_2(\cok(\phi')$ where $\I_j$ denotes the $j$th
Fitting ideal of a matrix.

\item If $M=\cok(phi)$ and $N=\cok(\psi)$, then $M \oplus N =\cok
\begin{pmatrix}\phi & 0 \\ 0 & \psi \end{pmatrix}$.

\item We will also need the fact that if $M=U \oplus V$ is a
non-trivial decomposition of a $2$-generated module $M$ over a
local ring, then both $U$ and $V$ are cyclic.

\end{enumerate}
\end{rmk}

\begin{prop}\label{cyclic}
Let $R$ be a reduced one-dimensional local ring. Let $M_R$ be a
MCM $R$-module. Then $(0:_RM)=\cap \Ass(M)$. In particular, if $M$
is cyclic, then $M \cong R/I$ where $I$ is an intersection of
minimal primes of $R$.
\end{prop}

\begin{proof}
Let $r \in (0:_RM)$ and let $P \in \Ass(M)$. Then there is an
injection $R/P \hookrightarrow M$ and thus $r(R/P)=0$; that is, $r
\in P$. As $P\in\Ass(M)$ was arbitrary, $r \in \cap \Ass(M)$. Now
let $r \in \cap \Ass(M)$. We want to show that $r \in (0:_R M)$.
Since $M$ is torsion-free it is enough to show that $(rM)_Q=0$ for
all minimal primes $Q$ of $R$. Since $R$ is one-dimensional and
$\m \not\in \Ass(M)$, we know that $\Ass(M)=\Supp(M)$. If $Q
\not\in \Ass(M)$ then $M_Q=0$ and thus $(rM)_Q=0$. On the other
hand, if $Q \in \Ass(M)$ then $r \in Q$ and as $R_Q$ is a field,
$rR_Q=0$. Therefore $(rM)_Q=0$ and hence $r \in (0:M)$.
\end{proof}

We are now ready to determine the MCM $R$-modules when $R$ is a
ring of type ($A2_n$), ($D2_n$), or ($D3$).

Suppose that $R$ is a ring of type ($A2_n$), a domain. Then by
(\ref{A2}), $S$ is a ring of type ($A_{2n+1}$) with minimal primes
$P_1=(x-\xi y^{n+1})$ and $P_2=(x+\xi y^{n+1})$. We know, from
\cite[9.9]{Yos90} and Lemma \ref{close}, that there are $n$
non-cyclic indecomposable MCM $S$-modules given by $2 \times 2$
matrix factorizations of $x^2-\xi^2y^{2n+2}$. Consider the matrix
factorizations $(\phi_j, \psi_j)$ for $x^2-\xi^2y^{2n+2}$:
\begin{equation}\label{a2n}\phi_j=\begin{bmatrix}-\xi^2y^{2n+2-j} & x \\ x &
-y^j\end{bmatrix}\text{\hspace{.25in}}(1 \leq j \leq
2n+2)\end{equation}

$$\psi_j=\begin{bmatrix}y^j & x \\ x & \xi^2y^{2n+2-j}\end{bmatrix}\text{\hspace{.25in}}(1 \leq j \leq 2n+2)$$

Let $M_i=\cok(\phi_j)$ and $N=\cok(\psi_j)$. It is is easy to see
that $M_j \cong N_j$, $M_j \cong M_{2n+2-j}$, and that $M_{n+1}
\cong S/(x-\xi y^{n+1}) \oplus S/(x+\xi y^{n+1})$. Since
$\I(M_j)=(x,y^j)$ ($1 \leq j \leq n$) is not equal to $\I(S/P_1)$,
$\I(S/P_2)$, or $\I(S/P_1)+I(S/P_2)$ we see by Remark (3) that
$M_j$ is indecomposable for each $j$, $1 \leq j \leq n$ and that
$M_i \not \cong M_j$ for $1 \leq i < j \leq n$. We note that each
$M_j$ is extended from an $R$-module since the entries of $\phi_j$
are monomials in generators for $\m$ with coefficients in $k$.
From Proposition \ref{galois}, $S/P_1$ and $S/P_2$ are not
extended and $\mu_M(S/P_1)=\mu_M(S/P_2)$ for all extended
$S$-modules $M$. Thus the indecomposable MCM $R$-modules extend to
the following $S$-modules:
\begin{equation}\label{a2modules}\text{$S$, $M_1$, $\ldots$ , $M_n$,
and $M_{n+1}\cong S/P_1 \oplus S/P_2$}\end{equation}

Now suppose that $R$ is a ring of type ($D2_n$). Then $R$ has two
minimal primes, $Q_1=(x-y)$ and $Q_2=(x^2-\xi^2y^{2n})$. From
(\ref{D2}), $S$ is a ring of type ($D_{2n+2}$) with three minimal
primes $P_1=(x-y)$, $P_2=(x-\xi y^n)$, and $P_3=(x+\xi y^n)$
satisfying $P_2\cap R=P_3\cap R=Q_2$. We now give matrix
factorizations of $\xi^2y^{2n}x-\xi^2y^{2n+1}-x^3+x^2y$ which will
in turn give us indecomposable MCM $S$-modules. Let $j$ run from
$1$ to $n$ and let $i$ run from $1$ to $n-1$. Now consider the
following $4n-2$ matrix factorizations $(\alpha_j, \beta_j)$,
$(\beta_j, \alpha_j)$, $(\phi_i, \psi_i)$, and $(\psi_i, \phi_i)$
where
\begin{eqnarray}\label{d2n}\alpha_j=\begin{bmatrix}\xi^2y^{2n+1-j} & x(x-y) \\ x &
y^{j-1}(x-y)\end{bmatrix}&& \beta_j=\begin{bmatrix}y^{j-1}(x-y) & -x(x-y) \\
-x & \xi^2y^{2n+1-j}\end{bmatrix}\end{eqnarray}
\begin{eqnarray*}\phi_i=\begin{bmatrix}\xi^2y^{2n-i} & x \\ x &
y^i\end{bmatrix}&& \psi_i=\begin{bmatrix}y^i(x-y) & -x(x-y)
\\ -x(x-y) & \xi^2y^{2n-i}(x-y)\end{bmatrix}\end{eqnarray*}

Now let $X_j=\cok(\alpha_j)$, $Y_j = \cok(\beta_j)$, $M_i
=\cok(\phi_i)$, and $N_i=\cok(\psi_i)$. We know from
\cite[9.12]{Yos90} and Lemma \ref{close} that $S$ has $4n-2$
non-cyclic indecomposable MCM modules and thus we need to justify
that the $4n-2$ modules given above are distinct and
indecomposable.

We note that $\I_2(\alpha_j)=\I_2(\beta_j)=(\xi^2y^{2n}-x^2)(x-y)$
for all $j$. Also, $\I_2(\phi_i)=(\xi^2y^{2n}-x^2)$ and
$\I_2(\psi_i)=(\xi^2y^{2n}-x^2)(x-y)^2$ for all $i$. Then by
Remark (3), for any $i\in \{1, \ldots , n-1\}$ and $j \in \{1,
\ldots , n\}$ we have that $X_i\not\cong M_j$, $X_i \not\cong
N_j$, and $M_i \not\cong N_j$. From Remark (2), $Y_j \cong
\syz_S^1(X_j)$. Then from \cite[9.12]{Yos90} we know that $X_j
\not\cong Y_j$ for all $j$.

Note that $\I_1(\alpha_j)=\I_1(\beta_j)=\I_1(\phi_j)=(x,y^j)$ and
that $\I_1(\psi_j)=(xy-x^2,xy^j-j^{j+1})$. Then by Remark (3) we
can conclude that $X_i \not \cong X_j$, $Y_i \not \cong Y_j$, $M_i
\not \cong M_j$, and $N_i \not \cong N_j$ for all $i \not= j$.
Thus the $4n-2$ MCM $S$-modules given by the matrix factorizations
in (\ref{d2n}) are distinct and we need only show they are
indecomposable.

The cyclic MCM $S$-modules are given by the cokernels of the
following six $1 \times 1$ matrices:
\begin{longtable}{lll}\label{cyclicideals}$A=(x-y)$ & $B=(\xi y^n-x)$ & $C=(\xi y^n+x)$ \\
$D=\left[(x-y)(\xi y^n-x)\right]$ & $E=\left[(x-y)(\xi
y^n+x)\right]$ & $F=(\xi^2y^{2n}-x^2)$\end{longtable}

Recall that $\I_1(X_j)=\I_1(Y_j)=(x,y^j)$ and
$\I_2(X_j)=\I_2(Y_j)=(\xi^2 y^{2n}-x^2)$. Also,
$\I_1(\phi_j)=(x,y^j)$, $\I_1(\psi_j)=(xy-x^2,xy^j-j^{j+1})$,
$\I_2(\phi_i)=(\xi^2y^{2n}-x^2)$ and
$\I_2(\psi_i)=(\xi^2y^{2n}-x^2)(x-y)^2$.

We note that
$\I_1(B)+\I_1(C)=(x,y^n)=\I_1(\alpha_n)=\I_1(\beta_n)$. However,
$$\I_2\left(\begin{bmatrix}B & 0 \\ 0 & C\end{bmatrix}\right) \not=
\I_2(\alpha_n)=\I_2(\beta_n)$$ and hence neither $X_n$ nor $Y_n$
is isomorphic to $\cok(B) \oplus \cok(C)$. Now from Remarks (3)
and (4), showing that the $4n-2$ MCM modules described in
(\ref{d2n}) are indecomposable reduces to showing that $(x,y^j)$
is not equal to the sum of any two of $A$, $B$, $C$, $D$, $E$,
$F$, not including $B+C$. It is easy to check that
$A+C+D+E+F=(x-y,\xi y^n+x)$ and $A+B+D+E+F=(x-y,\xi y^n-x)$,
neither of which contains $x$. Hence $X_j$, $Y_j$, $M_i$, and
$N_i$ are all indecomposable for $1 \leq i < j \leq n$.

Note that each of these modules is extended from a MCM $R$-module
since the entries of the matrices $\phi_j$, $\psi_j$, $\alpha_j$,
and $\beta_j$ are contained in $\m$ with coefficients in $k$.
Applying Proposition \ref{galois}, we see that the indecomposable
MCM $R$-modules extend to the following $S$-modules:
\begin{equation}\label{d2modules}\text{$M_1$, $\ldots$, $M_{n-1}$,
$N_1$, $\ldots$, $N_{n-1}$, $X_1$, $\ldots$, $X_n$, $Y_1$,
$\ldots$, $Y_n$,}\end{equation}
$$\text{$S$, $S/P_1$, $S/(P_2 \cap P_3)$, $A\cong S/P_2 \oplus S/P_3$, and
$B\cong S/(P_1 \cap P_2) \oplus S/(P_1 \cap P_3)$}$$

Finally, suppose that $R$ has type ($D3$). Then by (\ref{D3}) $S$
is a ring of type ($D_4$). Then $S$ has two non-cyclic
indecomposable MCM modules given by $2 \times 2$ matrix
factorizations of $y^3+ay^2x+byx^2+cx^3$. Consider the matrix
factorization $(\phi, \psi)$: \begin{equation}\label{d3matrices}\phi=\begin{bmatrix}y^2 & x^2  \\
-by-cx&y+ax\end{bmatrix}\text{ and }\psi=\begin{bmatrix}y+ax &
-x^2 \\ by+cx & y^2\end{bmatrix}\end{equation} Then $X=\cok(\phi)$
and $Y=\cok(\psi)$ are the non-cyclic indecomposable MCM
$S$-modules. We note that $X$ and $Y$ are extended since the
entries of $\phi$ and $\psi$ are contained in $\m$ with
coefficients in $k$. Applying Proposition \ref{galois} we see that
the indecomposable MCM $R$-modules extend to the following
$S$-modules
\begin{equation}\label{d3modules}\text{$S$, $X$, $Y$, $A\cong
S/P_1 \oplus S/P_2\oplus S/P_3$, and $B\cong S/(P_1 \cap P_2)
\oplus S/(P_1 \cap P_3) \oplus S/(P_2 \cap P_3)$}\end{equation}

We note that if $\phi$ is a matrix with entries in $\m$, the
maximal ideal of $R$, and $\cok(S^n \xrightarrow{\phi} S^n)$ is an
indecomposable MCM $S$-module, $\cok(R^n \xrightarrow{\phi} R^n)$
is an indecomposable MCM $R$-module. Thus the matrices given in
(\ref{a2n}), (\ref{d2n}), and (\ref{d3matrices}) give
presentations for indecomposable MCM $R$-modules where $R$ is a
ring of type ($A2_n$), ($D2_n$), and ($D3$) respectively.

We have now determined the indecomposable MCM $R$-modules when $R$
is a ring from Table \ref{hypersurfaces}. We turn now to the
non-hypersurfaces. Recall that if $R'$ is a complete reduced
equicharacteristic one-dimensional local non-hypersurface with
perfect residue field of characteristic not $2$, $3$, or $5$, $R'$
is isomorphic to $\End_R(\m)$ where $(R,\m)$ is one of the
hypersurfaces from Table \ref{hypersurfaces}.

Note that if $M$ is a MCM $R'$-module, $M$ has a natural
$R$-module structure since $R \subset R'$. Then as $M$ is
torsion-free as an $R'$-module, $M$ must be torsion-free as an
$R$-module. We also point out that each $R$-isomorphism is an
$R'$-isomorphism and that if $M$ is indecomposable as a
$R$-module, then $M$ is indecomposable as an $R'$-module. The
following result gives us a complete list of indecomposable MCM
$R'$-modules.

\begin{prop}\label{mod}
Let $R$ be a one-dimensional Gorenstein ring and assume that
$\End_R(\m)$ is local.

\begin{enumerate}
\item \cite[Prop. 7.2]{B63} Let $M$ be a MCM $R$-module and assume
that $M$ has no direct summand isomorphic to $R$. Then $M$ is a
MCM $\End_R(\m)$-module.

\item \cite[Lemma 9.4]{Yos90} Every indecomposable MCM
$\End_R(\m)$-module is an indecomposable MCM $R$-module.
\end{enumerate}

Note that $R$ does not have an $\End_R(\m)$-module structure since
$R \subsetneq \End_R(\m)$.
\end{prop}

Thus we now know that the indecomposable MCM $\End_R(\m)$-modules
are exactly the non-free indecomposable MCM $R$-modules.

Now we have succeeded in classifying all of the indecomposable MCM
modules over the rings classified in Theorem \ref{classification}.
When $R$ is a hypersurface with algebraically closed residue
field, the modules are listed in \cite[Chap. 9]{Yos90}. For the
rings of type ($A_n$), ($D_n$), ($E_6$), ($E_7$), and ($E_8$)
where the residue field is not necessarily algebraically closed we
appeal to Lemma \ref{close} and the classification in
\cite{Yos90}. When $R$ is a ring of type ($A2_n$), ($D2_n$), or
($D3$) the indecomposable modules are those that extend to the
modules listed in (\ref{a2modules}), (\ref{d2modules}), and
(\ref{d3modules}). When $R$ is a non-hypersurface, Proposition
\ref{mod} and the classifications of the modules over the
hypersurfaces give us the complete list of indecomposable MCM
$R$-modules.

\section{Ranks of Indecomposables}

Throughout this section we assume that $(R,\m,k)$ is one of the
complete local rings from Theorem \ref{classification}. Since $R$
is a one-dimensional CM local ring of finite CM type, $R$ is
reduced (cf. \cite{CWW}) and hence $R_P$ is a field for all
minimal primes $P$ of $R$. Recall that if $P$ is a minimal prime
of $R$ and $M$ is an $R$-module, the \emph{rank} $\rank_P(M)$
\emph{of} $M$ \emph{at} $P$ is defined to be the vector space
dimension of $M_P$ over the field $R_P$. If $P_1, \ldots P_s$ are
the minimal primes of $R$, the \emph{rank} of $M$ is the $s$-tuple
$\rank(M)=(r_1, \ldots, r_s)$ where $r_i=\rank_{P_i}(M)$.

We begin by calculating the ranks of the indecomposable modules
over rings of types ($A_n$), ($D_n$), ($E_6$), ($E_7$), and
($E_8$) using the AR-sequences given in \cite{Yos90}. We
illustrate this by computing the ranks for a ring $R$ of type
($D_6$). The remaining calculations are carried out in a similar
fashion in \cite{NB05}. The results of these calculations are
summarized in Theorem \ref{ranks} below.

Let $R$ be a ring isomorphic to $k[[x,y]]/(x^2y-y^5)$. Then from
\cite[9.12]{Yos90} we have the following AR-sequences:
\begin{center}\begin{tabular}{lcl} \ses {R/P_1} {X_1} {R/(P_2 \cap P_3)} & \hspace{.5in} & \ses {R/(P_2 \cap P_3)} {Y_1} {R/P_1} \\ \ses {X_1} {R/(P_2 \cap P_3) \oplus N_1} {Y_1} & & \ses {Y_1} {R/P_1 \oplus M_1} {X_1} \\ \ses {M_1} {X_1 \oplus Y_2} {N_1} & & \ses {N_1} {Y_1 \oplus X_2} {M_1} \end{tabular}\end{center}
For our purposes we need only know that these are short exact
sequences and that the rank function is additive along such
sequences. It is easy to see that $\rank(R)=(1,1,1)$ and that:
\begin{center}\begin{tabular}{lcl}$\rank(R/P_1)=(1,0,0)$ &
\hspace{.5in} & $\rank(R/(P_2 \cap P_3))=(0,1,1)$ \\
$\rank(R/P_2)=(0,1,0)$ & \hspace{.5in} & $\rank(R/(P_1 \cap P_3))=(1,0,1)$ \\
$\rank(R/P_3)=(0,0,1)$ & \hspace{.5in} & $\rank(R/(P_1 \cap
P_2))=(1,1,0)$
\end{tabular}\end{center}
Using this information along with the additivity of the rank
function along the AR-sequences, we see that
$\rank(X_1)=\rank(X_2)=\rank(Y_1)=\rank(Y_2)=(1,1,1)$,
$\rank(M_1)=(0,1,1)$, and $\rank(N_1)=(2,1,1)$.

We now determine the ranks of the indecomposable modules described
in the previous section; that is, we assume that $R$ is a ring of
type ($A2_n$), ($D2_n$) or ($D3$). We first state and prove the
following lemma.

\begin{lem}\label{tensorrank}
Let $R \rightarrow S$ be an extension of one-dimensional reduced
local rings. Let $P$ be a minimal prime of $R$ and let $Q$ be a
minimal prime of $S$ lying over $P$. If $_RN$ is a finitely
generated torsion-free $R$-module, then
$$\rank_P(N)=\rank_Q(N \otimes_RS)$$
\end{lem}

\begin{proof}
Let $r=\rank_P(N)$ and $s=\rank_Q(N \otimes_RS)$. Then $N_P \cong
R_P^r$ and $(N \otimes_RS)_Q \cong S_Q^s$. Now $(N\otimes_RS)_Q
\cong N_P \otimes_{R_P}S_Q \cong R_P^r \otimes_{R_P}S_Q \cong
S_Q^r$ and hence $S_Q^r \cong S_Q^s$. Thus $r=s$.
\end{proof}

Let $(R,\m,k)$ be a ring of type ($A2_n$). Then $S:= R \otimes_kK$
is a ring of type ($A_{2n+1}$). Then we know that the
indecomposable $R$-modules extend to the $S$-modules $S$, $M_1$,
$\ldots$ , $M_n$, and $M_{n+1}\cong S/P_1 \oplus S/P_2$, each
having rank $(1,1)$ over $S$. Since $R$ is a domain, both minimal
primes of $S$ lie over $(0)_R$. By Lemma \ref{tensorrank} we know
that each indecomposable MCM $R$-module has rank one.

Now suppose that $(R,\m,k)$ is a ring of type ($D2_n$). Then $S:=
R \otimes_kK$ is a ring of type ($D_{2n+2})$. Then we know that
the indecomposable $R$-modules extend to the $S$-modules $M_1$,
$\ldots$, $M_{n-1}$, $N_1$, $\ldots$, $N_{n-1}$, $X_1$, $\ldots$,
$X_n$, $Y_1$, $\ldots$, $Y_n$, $S$, $S/P_1$, $S/(P_2 \cap P_3)$,
$A\cong S/P_2 \oplus S/P_3$, and $B\cong S/(P_1 \cap P_2) \oplus
S/(P_1 \cap P_3)$. The ranks of these $S$-modules are $(0,1,1)$,
$\ldots$, $(0,1,1)$, $(2,1,1)$, $\ldots$, $(2,1,1)$, $(1,1,1)$,
$\ldots$, $(1,1,1)$, $(1,1,1)$, $\ldots$, $(1,1,1)$, $(1,1,1)$,
$(1,0,0)$, $(0,1,1)$, $(0,1,1)$, and $(2,1,1)$, respectively. The
ordering of the minimal primes is such that the second and third
lie over a common prime of $R$. Applying Lemma \ref{tensorrank},
we see that the ranks of the corresponding indecomposable MCM
$R$-modules are $(0,1)$, $\ldots$, $(0,1)$, $(2,1)$, $\ldots$,
$(2,1)$, $(1,1)$, $\ldots$, $(1,1)$, $(1,1)$, $\ldots$, $(1,1)$,
$(1,1)$, $(1,0)$, $(0,1)$, $(0,1)$, and $(2,1)$, respectively.

Finally, suppose that $R$ is a ring of type ($D3$), so that $S:=
R\otimes_kL$ is a ring of type ($D_4$). Then we know from the
previous section that the indecomposable MCM $R$-module extend to
the $S$-modules $S$, $X$, $Y$, $A\cong S/P_1 \oplus S/P_2\oplus
S/P_3$, and $B\cong S/(P_1 \cap P_2) \oplus S/(P_1 \cap P_3)
\oplus S/(P_2 \cap P_3)$. As $S$-modules, $S$, $X$, $Y$, and $A$
have constant rank $1$ and $B$ has constant rank $2$. Lemma
\ref{tensorrank} then gives us that the indecomposable MCM
$R$-modules have rank $1$, $1$, $1$, $1$, and $2$.

Now suppose that $R'$ is one of the non-hypersurfaces of finite CM
type. Then $R' \cong \End_R(\m)$ where $(R,\m)$ is a hypersurface
of finite CM type. By Proposition \ref{mod} the indecomposable MCM
$R'$-modules are exactly the non-free indecomposable MCM
$R$-modules. Thus we have already computed the ranks of the
indecomposable modules over the non-hypersurfaces. The following
theorem summarizes these calculations.

\begin{thm}\label{ranks} Let $R$ be one of the complete
hypersurfaces listed in Table \ref{hypersurfaces}. The ranks of
all indecomposable MCM modules over these rings are given in the
following table. In addition we give the number of indecomposables
of each rank.

\begin{center} {\small
\begin{tabular}{|ccc|ccc|ccc|}
\hline
\multicolumn{3}{|c}{one minimal prime} & \multicolumn{3}{c}{two minimal primes} & \multicolumn{3}{c|}{three minimal primes}\\
\hline  type & rank & \# & type & rank & \# & type & rank & \# \\
\hline $(A_{2n})$ & $1$ & $n+1$ & $(A_{2n+1})$ & $(1,0)$ & $1$ &
$(D_{2n+2})$ & $(1,0,0)$ & $1$ \\ \cline{1-3} $(E_6)$ & $1$ & $5$
& & $(0,1)$ & $1$ & & $(0,1,0)$ & $1$ \\ & $2$ & $2$ & & $(1,1)$ &
$n+1$ & & $(0,0,1)$ & $1$ \\ \cline{1-6} $(E_8)$ & $1$ & $7$ &
$(D_{2n+3})$ & $(1,0)$ & $1$ & & $(1,1,0)$ & $1$ \\ & $2$ & $7$ &
& $(0,1)$ & $n+1$ & & $(1,0,1)$ & $1$ \\ & $3$ & $3$ & & $(1,1)$ &
$2n+2$ & & $(0,1,1)$ & $n$ \\ \cline{1-3}
$(A2_n)$ & $1$ & $n+2$ & & $(2,1)$ & $n$ & & $(1,1,1)$ & $2n+1$ \\
\cline{1-6} $(D3)$ & $1$ & $4$ & $(E_7)$ & $(1,0)$ & $1$ & &
$(2,1,1)$ & $n-1$
\\ \cline{7-9} & $2$ & $1$ & & $(0,1)$ & $2$ & & & \\ \cline{1-3} & & & & $(1,1)$ & $6$ &
& &
\\ & & & & $(1,2)$ & $1$ & & & \\ & & & & $(2,1)$ & $2$ & & & \\ &
& & & $(2,2)$ & $3$ & & & \\ \cline{4-6} & & &
$(D2_n)$ & $(1,0)$ & $1$ & & & \\ & & & & $(0,1)$ & $n+1$ & & & \\
& & & & $(1,1)$ & $2n+1$ & & & \\ & & & & $(2,1)$ & $n$ & & & \\
\hline
\end{tabular}
}
\end{center}
\end{thm}

The reader may be disturbed by the lack of symmetry. For each of
the rings in Table \ref{hypersurfaces} we have fixed an order on
the minimal primes. This gives us, for example, an indecomposable
module of rank $(2,1)$ but no indecomposable module of rank
$(1,2)$ in the case when $R$ is a ring of type ($D_5$).

\section{MCM Modules over Non-complete Rings and a Krull-Schmidt Theorem}

Throughout this section $(R,\m,k)$ is an equicharacteristic
one-dimensional CM local ring with $k$ perfect of characteristic
not $2$, $3$ or $5$. We assume throughout that $R$ has FCMT,
equivalently \cite{RW94}, the $\m$-adic completion $\hat R$ has
FCMT. Thus $\hat R$ is isomorphic to a ring listed in Theorem
\ref{classification}. Since we do not have explicit formulas for
the rings $R$ we cannot give concrete descriptions of the modules
as in the complete case. Instead we represent $\C(R)$ as a
submonoid of the monoid $\C(\hat R)$ via the map $M \mapsto \hat
M$. In this paper we consider a \emph{monoid} to be a commutative,
cancellative, additive semigroup with $0$. We further restrict out
attention to reduced monoids --- monoids in which $0$ is the only
invertible element. Recall that $\bbn=\{0,1,2,\ldots\}$.

We put a monoid structure on $\C(R)$ by declaring
$[N]+[M]=[N\oplus M]$. As direct sum cancellation holds for
$R$-modules (\cite{E73}), $\C(R)$ satisfies our definition of a
monoid. As $\hat{R}$ is a complete local ring the decomposition of
$\hat{R}$-modules is unique up to isomorphism, \cite{A50}, and
hence $\C(\hat R)\cong \bbn^t$ where $t$ is the number of
isomorphism classes of indecomposable MCM $\hat{R}$-modules.  It
is shown in \cite{RW98} that the natural map taking $M$ to $\hat R
\otimes_R M$ induces a \emph{divisor homomorphism} $\C(R) \to
\C(\hat R)$, that is, for any two MCM $R$-modules $M$ and $N$, if
$\hat M\mid \hat N$ then $M\mid N$. Thus we may consider $\C(R)$
as a \emph{full submonoid} of $\bbn^t$; that is, a submonoid that
satisfies, for any $a,b \in \C(R)$, if $b=a+c$ for some $c \in
\bbn^t$ then $c \in \C(R)$.

Since we have a list (Theorem \ref{ranks}) of all of the
indecomposable MCM $\hat R$-modules as well as their ranks at each
of the minimal primes of $\hat R$, we can determine $\C(R)$ using
the following result, which follows as an immediate corollary to
\cite[Thm. 6.2]{LO96}:

\begin{prop}\label{LO}
Let $R$ and $\hat{R}$ be as above. In particular, $R$ is
one-dimensional and $\hat{R}$ is reduced. Let $M$ be a finitely
generated $\hat{R}$-module. Then $M$ is extended from an
$R$-module if and only if $\rank_P(M)=\rank_Q(M)$ whenever $P$ and
$Q$ are minimal primes of $\hat{R}$ lying over the same prime of
$R$.
\end{prop}

We will apply Proposition \ref{LO} to determine which of the MCM
$\hat{R}$-modules are extended from  MCM $R$-modules and then use
this information to determine $\C(R)$ as a full submonoid of
$\C(\hat{R})$. In order to determine the structure of $\C(R)$ we
need to introduce some additional terminology.

If there exists a divisor homomorphism $H \hookrightarrow \bbn^t$
we say that $H$ is a \emph{Krull monoid}. If $H \hookrightarrow
\bbn^t$ is a divisor homomorphism and each element of $\bbn^t$ is
the greatest lower bound of a finite set of elements of $\phi(H)$
then we say that $H \hookrightarrow \bbn^t$ is a \emph{divisor
theory} for $H$. In this case, we can define the \emph{class
group} $\Cl(H)$ of $H$ to be the cokernel of the induced map
$\Q(H) \hookrightarrow \bbz^t$, where $\Q(H)$ is the group of
formal differences of elements of $H$.

Let $H$ be a monoid and let $h \in H$. Then
$$L(h)=\{n\mid h=a_1+a_2+\cdots+a_n \text{ for irreducible $a_i$}\}$$ is the
\emph{set of lengths for the element} $h$ and $\mcL(H)=\{L(h)\mid
0 \not=h \in H\}$ is \emph{the system of lengths of} $H$. A monoid
$H$ is said to be \emph{factorial} if each element can be written
uniquely (up to order of the terms) as a sum of irreducible
elements of $H$. We note that this occurs exactly when $H$ is
free. A monoid $H$ is said to be \emph{half-factorial} if $L(h)$
is a singleton for each $h \in H$. The \emph{elasticity of an
element} $h \in H$ is $\rho(h)=\frac{\sup\{L(h)\}
}{\inf\{L(h)\}}$. The \emph{elasticity of the monoid} $H$ is
$\rho(H)=\sup\{\rho(h) \mid h \in H-\{0\}\}$. We note that
$\rho(H)=1$ if and only if $H$ is half-factorial.

The \emph{prime divisor classes} $G_0$ in $G=\Cl(H)$ are the
elements $q \in\Cl(H)$ such that $q=p+\Q(H)$ for some irreducible
element $p \in \bbn^t$. Note that the irreducible elements of
$\bbn^t$ are just the unit vectors $e_j$, $j=1, \ldots , t$.
Consider the map
$$\begin{array}{cccc}
  c: & \mcF(G_0) & \longrightarrow & G \\
     & \prod_{g \in G_0} g^{n_g} & \mapsto & \sum_{g \in G_0}n_gg \\
\end{array}$$ where $\mcF(G_0)$ is the free abelian monoid
(written multiplicatively) on the set $G_0$. The submonoid
$\mcB(G_0)=\{s \in \mcF(G_0):c(s)=0\} \subseteq \mcF(G_0)$ is
called the \emph{block monoid} of $G_0$. It is shown in \cite{G88}
that the system of lengths of $H$ is the same as the system of
lengths of $\mcB(G_0)$.

We now state and prove a lemma which allow us to calculate
$\Cl(\C(R))$ as well as the system of lengths of $\C(R)$.

\begin{lem}\label{calc}
\begin{enumerate}
\item \label{free} The divisor class group $\Cl(H)$ of a finitely
generated reduced Krull monoid $H$ is trivial if and only if
$H\cong \bbn^t$ for some $t$; i.e., $H$ is free.

\item \label{zero} Let $R$ and $\hat{R}$ be as above. If
$\#\Spec(\hat R)=\#\Spec(R)$, then $\C(R) \cong \C(\hat{R})$, and
hence $\Cl(\C(R))=0$.

\item \label{basis} If a monoid $H$ contains a $\bbz$-basis for a
group $G$ then $G=\Q(H)$.

\item Let $H=\ker(\mcA)\cap \bbn^t \subseteq \bbn^t$ where $\mcA$
is a $t \times s$ matrix with entries in $\bbz$. Further assume
that $H$ contains a $\bbz$-basis for $\ker(\mcA)$ and that the
natural inclusion $i:H \rightarrow \bbn^t$ is a divisor theory.
Then $\Cl(H)$ is isomorphic to the image of $\mcA:\bbz^s
\rightarrow \bbz^t$. Furthermore, the prime divisor classes in
$\Cl(H)$ are the elements $\{\mcA e_j\}_{j=1}^t$.
\end{enumerate}
\end{lem}

\begin{proof}
\begin{enumerate}

\item Suppose first that $\Cl(H)=0$. Then there is a divisor
theory $\phi:H \hookrightarrow \bbn^k$ with $\Q(H)=\bbz^k$.
Because $\phi$ is a divisor homomorphism we have that $$\phi(H)=
\Q(\phi(H))\cap\bbn^k=\bbz^k \cap \bbn^k=\bbn^k.$$

Suppose now that $H$ is free. Then there exists an isomorphism
$\phi:H \rightarrow \bbn^k$. Clearly $\phi$ is a divisor theory
and we have that $\Q(\phi):\Q(H) \rightarrow \bbz^k$ is also an
isomorphism. Thus $\Cl(H)=\cok(\Q(\phi))=0$.

\item Suppose that $\#\Spec(\hat R)=\#\Spec(R)$. Then each minimal
prime of $\hat R$ lies over a unique minimal prime of $R$. Thus it
is clear from Proposition \ref{LO} that all finitely generated
$\hat R$-modules are extended from $R$-modules. Thus $\C(R)
\hookrightarrow \C(\hat R)$ is an isomorphism. Then, as $\C(\hat
R)$ is free, so is $\C(R)$. By \ref{free} we have that
$\Cl(\C(R))=0$.

\item It is clear that $\Q(H) \subseteq G$. We now show the
reverse inclusion. If $\{h_1, \ldots, h_s\}$ is a $\bbz$-basis for
$G$ contained in $H$, then given $g \in G$ we can write
$g=a_1h_1+\cdots +a_sh_s$ with $a_i \in \bbz$. Re-index so that
for some $t\leq s$ we have $a_i \geq 0$ for $0 \leq i \leq t$ and
$a_i <0$ for $i$, $t < i \leq s$. Now we can write $g=(a_1h_1
+\cdots +a_th_t)-((-a_{t+1}h_{t+1})+\cdots +(-a_sh_s))$ making it
clear that $G \subseteq \Q(H)$.

\item By \ref{basis} we have that $\Q(H)=\ker(\mcA)$. Thus $\Cl(H)
\cong \bbz^t/\ker(\mcA)\cong \image(\mcA)$. It is now clear that
$\{\mcA e_j\}_{j=1}^t$ is the set of prime divisor classes in
$\Cl(H)$.
\end{enumerate}
\end{proof}

We now use the the lemma above to calculate $\C(R)$ and
$\Cl(\C(R))$ when $R$ is a non-complete local ring of FCMT. Here
we provide the calculations when $R$ is a local domain whose
completion is either a ring of type ($A_{2n+1}$) or a ring of type
($D_{2n+2}$). The remaining calculations are similar and are
worked out in detail in \cite{NB05}. Before we begin we recall the
following result, \cite[Thm. 1]{CL86}, of C. Lech.

\begin{prop}
Let $S$ be a complete Noetherian local ring. Then $S$ is the
completion of a Noetherian local domain if and only if the
following conditions hold.
\begin{enumerate}
\item The prime ring $\pi$ of $S$ is a domain and $S$ is a
torsion-free $\pi$-module. \item Either $S$ is a field or
$\depth(S) \geq 1$.
\end{enumerate}
\end{prop}

Since we are dealing only with equicharacteristic CM local rings
of dimension one these conditions are automatically satisfied.
Thus there is merit in the calculations that follow.

\subsubsection*{($A_{2n+1}$)}

Let $R$ be a local domain whose completion is isomorphic to the
ring $k[[x,y]]/(x^2-y^{2n})$. Referring to Proposition \ref{ranks}
we see that $\hat R$ has one indecomposable MCM module each of
rank $(0,1)$ and $(0,1)$ and $n+1$ indecomposable MCM modules of
rank $(1,1)$. Let $A$ be the indecomposable of rank $(0,1)$, let
$B$ the indecomposable of rank $(1,0)$ and let $M_0, \ldots, M_n$
be the indecomposables of constant rank one. If $L$ is any MCM
$\hat{R}$-module
$$L=\left(\bigoplus_{i=1}^n M_i^{m_i}\right) \oplus A^a \oplus B^b$$
where $m_i$, $a$, and $b$ are nonnegative integers and
$$\rank(L)=\left(\sum_{j=0}^n m_j+a,\sum_{j=0}^n m_j+b \right)$$

By Proposition \ref{LO}, $L$ is extended from a MCM $R$-module if
and only if $a=b$. Therefore the monoid of MCM $R$-modules is
$\C(R) \cong \bbn^{n+2}$. As $\C(R)$ is free, $\Cl(\C(R))=0$ by
Lemma \ref{calc}.

\subsubsection*{($D_{2n+2}$)}

Suppose now that $R$ is a local domain whose completion is
isomorphic to $k[[x,y]]/(x^2y-y^{2n+1})$. Referring to Proposition
\ref{ranks} we see that we have the following indecomposable MCM
$\hat R$-modules:

\begin{center}
\begin{tabular}{|cc|ccl|}
\hline module & rank & module & rank & \\ \hline
$A$ & $(1,0,0)$ & $E$ & $(1,0,1)$ & \\
$B$ & $(0,1,0)$ & $F_j$ & $(0,1,1)$ & $1 \leq j \leq n$ \\
$C$ & $(0,0,1)$ & $G_j$ & $(2,1,1)$ & $1 \leq j \leq n-1$ \\
$D$ & $(1,1,0)$ & $H_j$ & $(1,1,1)$ & $1 \leq j \leq 2n+1$ \\
\hline
\end{tabular}
\end{center}

Now if $L$ is any MCM $\hat R$-module,
$$L = A^a \oplus B^b \oplus C^c \oplus D^d \oplus E^e \oplus \left(\bigoplus_{j=1}^{n} F_j^{f_j}\right)\oplus \left(\bigoplus_{j=1}^{n-1} G_j^{g_j}\right)\oplus \left(\bigoplus_{j=1}^{2n+1} H_j^{h_j}\right)$$
Since $R$ is a domain, $L$ is extended if and only if $\rank(L)$
is constant; i.e., if
$$\sum_{j=1}^{n-1}g_j+a+e=\sum_{j=1}^{n}f_j+b\text{\ \  and\ \  }b+d=c+e.$$ Thus we have that
$\C(R)=\ker(\mcA)\cap \bbn^{2n+4} \oplus \bbn^{2n+1}$ where
$$\mcA=\begin{bmatrix}1&-1&0&0&1&-1&-1&\cdots &-1&1&\cdots &1\\0&1&-1&1&-1&0&0&\cdots &0&0&\cdots &0\end{bmatrix}$$

\vspace{.1in}

I claim that the following $2n+3$ elements form a $\bbz$-basis for
$\ker(\mcA)\oplus \bbz^{2n+1}$ that is contained in $\C(R)$.

\begin{eqnarray*}
\mathcal B= \left\{e_1+e_6,e_2+e_5,e_3+e_4,e_4+e_5+e_6\right\} &
\bigcup & \left\{e_1+e_j\mid 7\leq j \leq n+5\right\} \\ & \bigcup
& \left\{e_6+e_j \mid n+6 \leq 2n+4\right\} \\ & \bigcup &
\{e_j\mid 2n+5 \leq j \leq 4n+5\}
\end{eqnarray*}

Clearly $\mathcal B \subseteq \ker(\mcA)\oplus \bbz^{2n+1}$. Now
if $h=\sum_{i=1}^{4n+5}h_ie_i \in \ker(\mcA)$ then
$h_2+h_4=h_3+h_5$ and
$$h_1+h_5+\sum_{i=n+5}^{2n+4}h_i=h_2+h_6+\sum_{i=7}^{n+5}h_i$$ Thus
we can write $h$ uniquely as
\begin{eqnarray*}h & = &
\sum_{i=7}^{n+5}h_i(e_1+e_i)+\sum_{i=n+6}^{2n+4}h_i(e_6+e_i)
\\ & + & (h_1-\sum_{i=7}^{n+5}h_i)(e_1+e_6) \\ & + & (h_6-h_1+\sum_{i=7}^{n+5}h_i-\sum_{i=n+6}^{2n+4}h_i)(e_4+e_5+e_6) \\ & + & (h_4-h_6+h_1-\sum_{i=7}^{n+5}h_i+\sum_{i=n+6}^{2n+4}h_i)(e_3+e_4) \\ & + &
h_2(e_2+e_5)+\sum_{i=2n+5}^{4n+5}h_i(e_i)\end{eqnarray*}

Thus $\mathcal B$ is a $\bbz$-basis for $\ker(\mcA)\oplus
\bbz^{2n+1}$. One can check that each irreducible element of
$\bbn^{4n+5}$ is the greatest lower bound of two or three elements
of $\C(R)$. Since $e_j \in \C(R)$ for $2n+5 \leq j \leq 4n+5$, we
need only to consider $e_j$ for $j \leq 2n+4$. Let $j$ and $k$ be
such that $7 \leq j \leq n+5$, and $n+6 \leq k \leq 2n+4$. Then
$e_1=\glb(e_1+e_6,e_1+e_2+e_3), e_2=\glb(e_2+e_5,
e_1+e_2+e_3),e_3=\glb(e_3+e_4, e_1+e_2+e_3)$, $e_4=\glb(e_3+e_4,
e_4+e_5+e_6)$, $e_5=\glb(e_2+e_5,e_4+e_5+e_6), e_6=\glb(e_1+e_6,
e_4+e_5+e_6)$, $e_j=\glb(e_1+e_j, e_j+e_4+e_5)$, and
$e_k=\glb(e_6+e_k, e_k+e_2+e_3)$. Thus, the natural inclusion
$\C(R) \hookrightarrow \C(\hat{R})$ is a divisor theory. By Lemma
\ref{calc}, $\Cl(\C(R))\cong \bbz\oplus \bbz$.

The remaining calculations may be worked out in a similar fashion,
and the following proposition summarizes the results.

\begin{prop}\label{monoids} Let $R$ be a one-dimensional local ring whose
completion is equicharacteristic with perfect residue field of
characteristic different from $2$, $3$ and $5$. Then $\C(R)$ and
$\Cl(\C(R))$ depend only on the singularity type of $\hat{R}$ and
on $m:=\#\Spec(\hat{R})-\#\Spec(R)$. The results are summarized in
the following table. We list the results only when $m>0$ since the
case $m=0$ was taken care of in Lemma \ref{calc}.

\begin{center}
\begin{longtable}{|c|c|l|l|l|}\caption{Monoids of  MCM Modules}\\
\hline
 $\hat{R}$ & $m$ & $\C(\hat{R})$ & $\C(R)$ & $\Cl(\C(R))$ \\
\hline \endhead \hline \endfoot

$(A_{2n+1})$& $1$ & $\bbn^{n+3}$ & $\bbn^{n+2}$ & $0$
\\

$(D_{2n+3})$& $1$ & $\bbn^{4n+4}$ & $\ker(\mcA_1)\cap
\bbn^{2n+2}\oplus \bbn^{2n+2}$ & $\bbz$ \\

$(D'_{2n+3})$& $1$ & $\bbn^{4n+3}$ & $\ker(\mcA_1)\cap
\bbn^{2n+2}\oplus \bbn^{2n+1}$ & $\bbz$ \\

$(D_{2n+2})$$^1$& $1$ & $\bbn^{4n+5}$ & $\ker(\mcA_2)\cap
\bbn^{2n+2}\oplus \bbn^{2n+3}$ & $\bbz$ \\

$(D'_{2n+2})$$^1$& $1$ & $\bbn^{4n+4}$ & $\ker(\mcA_2)\cap
\bbn^{2n+2}\oplus \bbn^{2n+2}$ & $\bbz$ \\

$(D_{2n+2})$$^2$& $1$ & $\bbn^{4n+5}$ & $\ker(\mcA_3)\cap \bbn^4\oplus \bbn^{4n+1}$ & $\bbz$ \\

$(D'_{2n+2})$$^2$& $1$ & $\bbn^{4n+4}$ & $\ker(\mcA_3)\cap \bbn^4\oplus \bbn^{4n}$ & $\bbz$ \\

$(D_{2n+2})$& $2$ & $\bbn^{4n+5}$ & $\ker(\mcA_4)\cap
\bbn^{2n+4} \oplus \bbn^{2n+1}$ & $\bbz\oplus \bbz$ \\

$(D'_{2n+2})$& $2$ & $\bbn^{4n+4}$ & $\ker(\mcA_4)\cap
\bbn^{2n+4} \oplus \bbn^{2n}$ & $\bbz\oplus \bbz$ \\

$(E_7)$& $1$ & $\bbn^{15}$ & $\ker(\mcA_5) \cap \bbn^6 \oplus \bbn^9$ & $\bbz$ \\

$(E'_7)$& $1$ & $\bbn^{14}$ & $\ker(\mcA_5) \cap \bbn^6 \oplus \bbn^8$ & $\bbz$ \\

$(D2_n)$& $1$ & $\bbn^{4n+3}$ & $\ker(\mcA_6)\cap\bbn^{2n+2}\oplus \bbn^{2n+1}$ & $\bbz$ \\

$(D2'_n)$& $1$ & $\bbn^{4n+2}$ & $\ker(\mcA_6)\cap\bbn^{2n+2}\oplus \bbn^{2n}$ & $\bbz$ \\

\hline

\end{longtable}
\end{center}

The matrices in the table above are as follows:

\begin{itemize}

\item $\mcA_1=\begin{bmatrix}
1&-1&|&-1&\cdots&-1&1&\cdots&1\end{bmatrix}_{1 \times (2n+2)}$

\item
$\mcA_2=\begin{bmatrix}1&-1&1&0&-1&0&\vline&-1&\cdots&-1&1&\cdots&1\end{bmatrix}_{1
\times (2n+2)}$

\item $\mcA_3=\begin{bmatrix}-1&1&1&-1\end{bmatrix}_{1 \times 4}$

\item $\mcA_4=\begin{bmatrix}1&-1&0&0&1&-1&-1&\cdots &-1&1&\cdots
&1\\0&1&-1&1&-1&0&0&\cdots &0&0&\cdots &0\end{bmatrix}_{2 \times
(2n+4)}$

\item $\mcA_5=\begin{bmatrix}1&-1&1&-1&-1&1\end{bmatrix}$

\item $\mcA_6=\begin{bmatrix}
1&-1&|&-1&\cdots&-1&1&\cdots&1\end{bmatrix}_{1 \times (2n+2)}$

\end{itemize}

{\bf Note:} The two cases for $(D_{2n+2})$ and $(D'_{2n+2})$ and
$m=1$ correspond to the cases 1: $P_1 \cap R= P_2 \cap R$ and 2:
$P_2 \cap R = P_3 \cap R$. As expected, the two cases coincide
when $n=4$.

\end{prop}

We are now ready to prove our main result.

\begin{thm}\label{main}
Let $R$ be a reduced one-dimensional equicharacteristic local ring
with perfect residue field of characteristic not $2$, $3$, or $5$.
Suppose further that $R$ has finite Cohen-Macaulay type. Let
$m=\#\Spec(\hat{R})-\#\Spec(R)$.
\begin{enumerate}
\item If $m=0$ then $\C(R)$ is factorial and Krull-Schmidt holds
for torsion-free $R$-modules.

\item Let $m=1$. If $\hat{R}$ is a ring of type $(A_{2n+1})$ then
$\C(R)$ is factorial; otherwise $\C(R)$ is half-factorial but not
factorial.

\item If $m=2$ then $\C(R)$ is not half-factorial.
\end{enumerate}
\end{thm}

\begin{proof}

First recall from part \ref{zero} of Lemma \ref{calc} that
$\Cl(\C(R))=0$ if and only if $\C(R)$ is free and hence factorial.

We now turn to the case where $\Cl(\C(R)) \cong \bbz$. This is the
case when $R$ is a domain and $\hat R$ is a ring of type
($D_{2n+2}$), ($D_{2n+3}$), ($E_7$), or ($D2_n$) with
$\#\Spec(\hat R)-\#\Spec(R)=1$. From Proposition \ref{monoids} we
know that in each of these cases $\C(R) \cong \ker (\mcA) \cap
\bbn^s \oplus \bbn^t$ for some $1 \times s$ integer matrix $\mcA$
where $s$ and $t$ are non-negative integers.

From Lemma \ref{calc} we know that in order to determine the block
monoid of $\C(R)$ we need only compute $\mcA'(e_j)$ for all $j$,
$1 \leq j \leq s$, where $\mcA'=\begin{bmatrix}\mcA &\mid & 0 &
\cdots & 0\end{bmatrix}_{1 \times (s+t)}$. Referring to the
matrices $\mcA$ in Proposition \ref{monoids} we see that the set
of prime divisor classes of $\Cl(\C(R))$ is $G_0=\{0,+1,-1\}$ and
hence $\mcB(\C(R))\cong \bbn^2$. Since $\mcL(\bbn^2)=\{\{n\}\mid n
\in \bbn\}$ and as the system of lengths for $\C(R)$ is
necessarily the same as for $\mcB(\C(R))$ we have that $\C(R)$ is
half-factorial. We note that since $\Cl(\C(R)) \not=0$, $\C(R)$ is
not factorial by part \ref{free} of Lemma \ref{calc}.

We now compute the block monoid associated to $\Cl(\C(R)) \cong
\bbz \oplus \bbz$, which occurs when $R$ is a domain and $\hat{R}$
is of type ($D_{2n+2}$). Recall that $\C(R) \cong \ker(\mcA) \cap
\bbn^{2n+4}\oplus \bbn^{2n+1}$ where
\begin{equation}\label{matrix}\mcA=\begin{bmatrix}1&-1&0&0&1&-1&-1&\cdots
&-1&1&\cdots &1\\0&1&-1&1&-1&0&0&\cdots &0&0&\cdots
&0\end{bmatrix}.\end{equation}

By Lemma \ref{calc} the set of the prime divisor classes is
$$\{(1,0),(-1,0),(0,1),(0,-1),(1,-1),(-1,1)\}$$ and thus the block monoid
is

$$\mcB(\C(R))\cong \ker\begin{bmatrix}1&-1&0&0&1&-1 \\ 0&1&-1&1&-1&0\end{bmatrix}\bigcap \bbn^6.$$

The irreducible elements of this monoid are:
\begin{center}
\begin{tabular}{lll}
$h_1=(1,0,0,0,0,1)$ & $h_2=(0,1,0,0,1,0)$ & $h_3=(0,0,1,1,0,0)$ \\
$h_4=(1,1,1,0,0,0)$ & $h_5=(0,0,0,1,1,1)$ & \\
\end{tabular}
\end{center}

Note that $h_1+h_2+h_3=h_4+h_5$ and hence this monoid is not
half-factorial.

\end{proof}

We conclude this section with a result which shows how the class
of torsion-free modules behaves much better than the class of
arbitrary finitely generated modules. This is exhibited using the
elasticity function $\rho$.

\begin{thm}
Let $R$, $\hat R$, and $m$ be as in Theorem \ref{main}. Let
$\M(R)$ denote the monoid of all finitely generated $R$-modules.
\begin{enumerate}
\item \label{KS} If $m=0$ then the Krull-Schmidt property holds
for the class of all finitely generated $R$-modules.

\item \label{m=1} Let $m=1$. \begin{enumerate} \item If $\hat R$
is a ring of type ($A_{2n+1}$) for some $n\geq 0$, then $\C(R)$ is
factorial; otherwise $\rho(\C(R))=1$ but $\C(R)$ is not factorial.
\item If $\hat R$ is a ring of type ($A_1$) then $\rho(\M(R))=1$
but $\M(R)$ is not factorial; otherwise
$\rho(\M(R))=\infty$.\end{enumerate}

\item \label{m=2} If $m=2$, then $\rho(\C(R))=\frac{3}{2}$ and
$\rho(\M(R))=\infty$.
\end{enumerate}
\end{thm}

\begin{proof}
Since Proposition \ref{LO} holds for all finitely generated
modules, part \ref{zero} of Lemma \ref{calc} gives \ref{KS}.

We now suppose that $m>0$. A result in \cite{HKKW} says that if
$\hat R$ is a ring from Table \ref{hypersurfaces} not of type
($A_1$) (in other words, not Dedekind-like) then for any $s$-tuple
of natural numbers $(r_1, \ldots , r_s)$ (where $s$ is the number
of minimal primes of $\hat R$) there exists an indecomposable
finitely generated $\hat R$-module $M$ such that the torsion-free
rank of $M$ is $(r_1, \ldots , r_s)$.

Suppose now that $R$ is a domain and $\hat R$ has two minimal
primes. Let $n \in \bbn$ and set $N=1+2+\cdots +n$. From
\cite{HKKW} we have that for all $i \in \{1,2, \ldots , n, N\}$
there exist indecomposable finitely generated $\hat R$-modules
$M_{0,i}$ and $M_{i,0}$ of rank $(0,i)$ and $(i,0)$, respectively.
Since $R$ is a domain, none of these modules is extended. However,
$M_{0,1}\oplus M_{1,0}, \ldots, M_{0,n}\oplus M_{n,0}$,
$\left(\bigoplus_{i =1}^n M_{i,0}\right)\oplus M_{0,N}$, and
$\left(\bigoplus_{i =1}^n M_{0,i}\right)\oplus M_{N,0}$ are
extended. Moreover, none of these extended modules has a direct
summand which is extended. Therefore, they are extended from
indecomposable $R$-modules. Now since we have
\begin{eqnarray*} M & = & \left(\left(\bigoplus_{i =1}^n M_{i,0}\right)\oplus
M_{0,N}\right)\bigoplus \left(\left(\bigoplus_{i =1}^n
M_{0,i}\right)\oplus M_{N,0}\right) \\ & \cong & \bigoplus_{i=1}^n
\left(M_{0,i}\oplus M_{i,0}\right) \bigoplus \left(M_{0,N}\oplus
M_{N,0}\right)\end{eqnarray*} we see that $2, n+1 \in L([M])$ and
hence the elasticity of the element $[M] \in \M(R)$ is at least
$\frac{n+1}{2}$. As $n$ was arbitrary the elasticity of $\M(R)$ is
infinite. A similar argument takes care of the remaining cases
when $m>0$ and $\hat R$ is not of type ($A_1$).

Now suppose that $R$ is a local domain whose completion $\hat R$
is of type ($A_1$). Then $\hat R$ is split Dedekind-like and hence
has infinitely many non-isomorphic indecomposables finitely
generated modules of ranks $(1,0)$ and $(0,1)$. Moreover, all
indecomposable finitely generated $\hat R$-modules have ranks
$(1,0)$, $(0,1)$, or $(1,1)$ (c.f. \cite{KL2}). Thus there are
non-isomorphic finitely generated indecomposable modules $M$ and
$M'$ with rank $(1,0)$ and non-isomorphic indecomposable modules
$N$ and $N'$ with rank $(0,1)$. None of these modules is extended,
but $M\oplus N$, $M'\oplus N'$, $M \oplus N'$, and $M'\oplus N$
are extended. Thus we have
$$\left(M \oplus M'\right) \bigoplus \left(N \oplus N'\right)
\cong \left(M \oplus N'\right) \bigoplus \left(M' \oplus
N\right)$$ which exhibits the failure of Krull-Schmidt for
finitely generated $R$-modules. However, since each extended $\hat
R$-module $M$ must have constant rank, any direct sum
decomposition of $M$ must have the same number of indecomposable
summands of rank $(1,0)$ as it has of rank $(0,1)$. Therefore,
$\M(R)$ is half-factorial and $\rho(\M(R))=1$.

Now we deal with the monoid $\C(R)$. When $m\leq 1$ Theorem
\ref{main} gives us that $\C(R)$ is half-factorial and hence
$\rho(\C(R))=1$. Now suppose $m=2$. Recall that the set of lengths
of $\C(R)$ is the same as the set of lengths of $H \cong
\ker(\mcA) \cap \bbn^6$ where the matrix $\mcA$ is as in
(\ref{matrix}). Using the algorithm from \cite[Sec. 2]{K04} we
easily determine that $\rho(\C(R))=\rho(H)=\frac{3}{2}$.
\end{proof}

\section{No Indecomposable of Rank $4$}\label{norank4}

\indent In the previous section we computed the monoid of MCM
$R$-modules if $R$ is an equicharacteristic one-dimensional local
domain of finite representation type. From the descriptions of the
monoids given it is easy to see that the possible ranks of
indecomposable MCM $R$-modules over a domain $R$ are $1$, $2$, and
$3$. This result is in contradiction with a result in \cite{RSW94}
which states that the ring $R=\mathbb Q[[T,\sqrt[3]{2}T]]$ has an
indecomposable torsion free module of rank $4$. We note that this
ring is a ring of type ($D3'$). We first recall the basic set-up
in \cite{RSW94}.

Let $R$ be a reduced one-dimensional local ring with module-finite
integral closure $\bar R\not=R$. Let $\f$ be the conductor ideal.
Then $R$ can be represented as the pullback in the following
diagram:

\centerline{\xymatrix{R \ar[r] \ar@{>>}[d] & \bar R \ \ar@{>>}[d]
\\ R/\f \ar[r] & \bar R/\f \\}}

The bottom line of the pullback diagram is called an
\emph{Artinian pair}. (The rings are Artinian since $\f$ contains
a non-zero-divisor.) A module over an Artinian pair $A \rightarrow
B$ is a pair $V \rightarrow W$ such that $W$ is a finitely
generated projective $B$-module and $V$ is an $A$-submodule of $B$
satisfying $BV=W$. A direct sum decomposition of $V \rightarrow W$
is a direct sum decomposition of $W$ as a $B$-module, $W=W_1\oplus
\cdots \oplus W_n$, such that $V=(W_1\cap V) \oplus \cdots \oplus
(W_n \cap V)$.

When $R=\mathbb Q[[T,\sqrt[3]{2}T]]$ and $\bar R= \mathbb Q
(\sqrt[3]{2})[[T]]$ the Artinian pair is $\mathbb Q \rightarrow
\mathbb Q (\sqrt[3]{2})$. By Proposition 2.2 of \cite{RSW94} there
is a bijective correspondence between the indecomposable modules
over this Artinian pair and the indecomposable MCM $R$-modules.

We now show how the supposed indecomposable $M$ of rank $4$ in
\cite{RSW94} decomposes as a direct sum of two modules of rank
$2$. We have that $k=\q$, $K=\q(\sqrt[3]{2})$ and
$L=\q(\sqrt[3]{2},\omega)$ where $\omega$ is a primitive cube root
of unity. The Galois group of $L/k$ is generated by an element
$\sigma$ of order two and an element $\tau$ of order three. We
wish to decompose the module $N=(L \xrightarrow{\delta} L \times
L)$ ($\delta(x)=(x,x)$) as a $k \rightarrow K$ module where the
$K$-action on $L \times L$ is defined to be $s \cdot (x,y) =
(s^{\tau}x,sy)$ for all $s \in K$ and $(x,y) \in L \times L$. It
is easier to understand the decomposition when we instead consider
the isomorphic $k\rightarrow K$ module $N=(L \xrightarrow{\gamma}
L \times L)$, where now ($\gamma(x)=(x^{\tau^2},x$) and $K$ acts
diagonally on $L \times L$: $s \cdot (x,y) = (sx,sy)$ for all $s
\in K$ and $(x,y) \in L \times L$. This $k \rightarrow K$ module
isomorphism is given by the following commutative diagram where
the top row represents the module listed in \cite{RSW94} and
$\phi(x,y)=(x^{\tau^2},y)$.

\centerline{\xymatrix{M: & L \ar[r]^{\delta} \ar[d]_{id} & L
\times L \ar[d]^{\phi} \\ N: & L \ar[r]_{\gamma} & L \times L \\}}

The module $L \rightarrow L \times L$, as given by the bottom row
of the diagram, decomposes over $k \rightarrow K$ as
$$\{(x^{\tau^2}+y^{\tau^2}\omega^2,x+y\omega) \mid x,y \in K\}
\bigoplus \{(x^{\tau^2}\omega+y^{\tau^2},x\omega+y\omega^2)\mid
x,y \in K\}$$ The details are straightforward to check and can be
found in \cite{NB05}. We also note that the supposed
indecomposable of rank three over this ring, $K \rightarrow L
\times K$, decomposes over $k \rightarrow K$ as a direct sum of
three rank one modules, each isomorphic to the module $k
\rightarrow K$.

\bibliographystyle{amsplain}
\bibliography{references}

\end{document}